# VARIANCE ESTIMATION IN NONPARAMETRIC REGRESSION VIA THE DIFFERENCE SEQUENCE METHOD

By Lawrence D. Brown[1] and M. Levine[1,2]

*University of Pennsylvania and Purdue University*

Consider a Gaussian nonparametric regression problem having both an unknown mean function and unknown variance function. This article presents a class of difference-based kernel estimators for the variance function. Optimal convergence rates that are uniform over broad functional classes and bandwidths are fully characterized, and asymptotic normality is also established. We also show that for suitable asymptotic formulations our estimators achieve the minimax rate.

**1. Introduction.** Let us consider the nonparametric regression problem

$$(1) \qquad y_i = g(x_i) + \sqrt{V(x_i)}\epsilon_i, \qquad i=1,\ldots,n,$$

where $g(x)$ is an unknown mean function, the errors $\epsilon_i$ are i.i.d. with mean zero, variance 1 and the finite fourth moment $\mu_4 < \infty$ while the design is fixed. We assume that $\max\{x_{i+1} - x_i\} = O(n^{-1})$ for $\forall i = 0, \ldots, n$. Also, the usual convention $x_0 = 0$ and $x_{n+1} = 1$ applies. The problem we are interested in is estimating the variance $V(x)$ when the mean $g(x)$ is unknown. In other words, the mean $g(x)$ plays the role of a nuisance parameter. The problem of variance estimation in nonparametric regression was first seriously considered in the 1980s. The practical importance of this problem has been also amply illustrated. It is needed to construct a confidence band for any mean function estimate (see, e.g., Hart [24], Chapter 4). It is of interest in confidence interval determination for turbulence modeling (Ruppert et al. [34]), financial time series (Härdle and Tsybakov [23], Fan and Yao [18]), covariance structure estimation for nonstationary longitudinal data (see, e.g.,

Received August 2006; revised December 2006.
[1]Supported in part by NSF Grant DMS-04-05716.
[2]Supported in part by a 2004 Purdue Research Foundation Summer Faculty grant.

Diggle and Verbyla [10]), estimating correlation structure of heteroscedastic spatial data (Opsomer et al. [31]), nonparametric regression with lognormal errors as discussed in Brown et al. [2] and Shen and Brown [36], and many other problems.

In what follows we describe in greater detail the history of a particular approach to the problem. von Neumann [40, 41] and then Rice [33] considered the special, homoscedastic situation in which $V(x) \equiv \sigma^2$ in the model (1) but $\sigma^2$ is unknown. They proposed relatively simple estimators of the form

$$\hat{V}(x) = \frac{1}{2(n-1)} \sum_{i=1}^{n-1} (y_{i+1} - y_i)^2. \tag{2}$$

The next logical step was made in Gasser, Sroka and Jennen-Steinmetz [19], where three neighboring points were used to estimate the variance,

$$\hat{V}(x) = \frac{2}{3(n-2)} \sum_{i=1}^{n-2} \left(\frac{1}{2} y_i - y_{i+1} + \frac{1}{2} y_{i+2}\right)^2. \tag{3}$$

A further general step was made in Hall, Kay and Titterington [21]. The following definition is needed first.

DEFINITION 1.1. Let us consider a sequence of numbers $\{d_i\}_{i=0}^{r}$ such that

$$\sum_{i=0}^{r} d_i = 0 \tag{4}$$

while

$$\sum_{i=0}^{r} d_i^2 = 1. \tag{5}$$

Such a sequence is called a difference sequence of order $r$.

For example, when $r = 1$, we have $d_0 = \frac{1}{\sqrt{2}}, d_1 = -d_0$, which defines the first difference $\Delta Y = \frac{Y_i - Y_{i-1}}{\sqrt{2}}$. The estimator of Hall, Kay and Titterington [21] can be defined as

$$\hat{V}(x) = (n-r)^{-1} \sum_{i=1}^{n-r} \left(\sum_{j=0}^{r} d_j y_{j+i}\right)^2. \tag{6}$$

The conditions (4) and (5) are meant to insure the unbiasedness of the estimator (6) when $g$ is constant and also the identifiability of the sequence $\{d_i\}$.

A different direction was taken in Hall and Carroll [20] and Hall and Marron [22] where the variance was estimated by an average of squared residuals from a fit to $g$; for other work on constant variance estimation, see also Buckley, Eagleson and Silverman [5], Buckley and Eagleson [4] and Carter and Eagleson [7].

The difference sequence idea introduced by Hall, Kay and Titterington [21] can be modified for the case of a nonconstant variance function $V(x)$. As a rule, the average of squared differences of observations has to be localized in one way or another—for example, by using the nearest neighbor average, a spline approach or local polynomial regression. The first to try to generalize it in this way were probably Müller and Stadtmüller [27]. It was further developed in Hall, Kay and Titterington [21], Müller and Stadtmüller [28], Seifert, Gasser and Wolf [35], Dette, Munk and Wagner [9], and many others. An interesting application of this type of a variance function estimator for the purpose of testing the functional form of the given regression model is given in Dette [8].

Another possible route to estimating the variance function $V(x)$ is to use the local average of the squared residuals from the estimation of $g(x)$. One of the first applications of this principle was in Hall and Carroll [20]. A closely related estimator was also considered earlier in Carroll [6] and Matloff, Rose and Tai [26]. This approach has also been considered in Fan and Yao [18].

Some of the latest work in the area of variance estimation includes attempts to derive methods that are suitable for the case where $X \in \mathcal{R}^d$ for $d > 1$; see, for example, Spokoiny [38] for generalization of the residual-based method and Munk, Bissantz, Wagner and Freitag [29] for generalization of the difference-based method.

The present research describes a class of nonparametric variance estimators based on difference sequences and local polynomial estimation, and investigates their asymptotic behavior. Section 2 introduces the estimator class and investigates its asymptotic rates of convergence as well as the choice of the optimal bandwidth. Section 3 establishes the asymptotic normality of these estimators. Section 4 investigates the question of asymptotic minimaxity for our estimator class among all possible variance estimators for nonparametric regression.

**2. Variance function estimators.** Consider the model (1). We begin with the following formal definition.

DEFINITION 2.1. A pseudoresidual of order $r$ is

(7) $$\Delta_i \equiv \Delta_{r,i} = \sum_{j=0}^{r} d_j y_{j+i-\lfloor r/2 \rfloor},$$

where $\{d_j\}$ is a difference sequence satisfying (4)–(5) and $i = \lfloor \frac{r}{2} \rfloor + 1, \ldots, n + \lfloor \frac{r}{2} \rfloor - r$.

Let $K(\cdot)$ be a real-valued function such that $K(u) \geq 0$ and is not identically zero; $K(u)$ is bounded [$\exists M > 0$ such that $K(u) \leq M$ for $\forall u$]; $K(u)$ is supported on $[-1, 1]$ and $\int K(u) \, du = 1$. We use the notation $\sigma_K^2 = \int u^2 K(u) \, du$ and $R_K = \int K^2(u) \, du$. Then, based on $\Delta_{r,i}$, we define a variance estimator $\hat{V}_h(x)$ of order $r$ as the local polynomial regression estimator based on $\Delta_{r,i}^2$,

$$\hat{V}_h(x) = \hat{a}_0, \tag{8}$$

where

$$(\hat{a}_0, \hat{a}_1, \ldots, \hat{a}_p)$$
$$= \arg\min_{a_0, a_1, \ldots, a_p} \sum_{i=\lfloor r/2 \rfloor + 1}^{n + \lfloor r/2 \rfloor - r} [\Delta_{r,i}^2 - a_0 - a_1(x - x_i) - \cdots - a_p(x - x_i)^p]^2$$
$$\times K\left(\frac{x - x_i}{h}\right).$$

The value $h$ in (8) is called the bandwidth and $K$ is the weight function.

It should be clear that these estimators are unbiased under the assumption of homoscedasticity $V(x) \equiv \sigma^2$ and constant mean $g(x) \equiv \mu$. We begin with the definition of the functional class that will be used in the asymptotic results to follow.

DEFINITION 2.2. Define the functional class $\mathcal{C}_\gamma$ as follows. Let $C_1 > 0$, $C_2 > 0$. Let us denote $\gamma' = \gamma - \lfloor \gamma \rfloor$ where $\lfloor \gamma \rfloor$ denotes the greatest integer less than $\gamma$. We say that the function $f(x)$ belongs to the class $\mathcal{C}_\gamma$ if for all $x, y \in (0, 1)$

$$|f^{\lfloor \gamma \rfloor}(x) - f^{\lfloor \gamma \rfloor}(y)| \leq C_1 |x - y|^{\gamma'}, \tag{9}$$

$$|f^{(k)}(x)| \leq C_2, \tag{10}$$

for $k = 0, \ldots, \lfloor \gamma \rfloor - 1$. Note that $\mathcal{C}_\gamma$ depends on the choice of $C_1$, $C_2$, but for our convenience we omit this dependence from the notation. There are also similar types of dependence in the definitions that immediately follow.

DEFINITION 2.3. Let $\delta > 0$. We say the function is in class $\mathcal{C}_\gamma^+$ if it is in $\mathcal{C}_\gamma$ and in addition

$$f(x) \geq \delta. \tag{11}$$

These classes of functions are familiar in the literature, as in Fan [15, 16] and are often referred to as Lipschitz balls.

DEFINITION 2.4. Define the pointwise risk of the variance estimator $\hat{V}_h(x)$ (its mean squared error at a point $x$) as

$$R(V(x), \hat{V}_h(x)) = E[\hat{V}_h(x) - V(x)]^2.$$

DEFINITION 2.5. Define the global mean squared risk of the variance estimator $\hat{V}_h(x)$ as

(12) $$R(V, \hat{V}_h) = E\left(\int_0^1 (\hat{V}_h(x) - V(x))^2 \, dx\right).$$

Then the globally optimal in the minimax sense bandwidth $h_{\text{opt}}$ is defined as

$$h_n = \arg\min\{\sup\{R(V, \hat{V}_h) : V \in \mathcal{C}_\gamma, g \in \mathcal{C}_\beta\} : h > 0\}.$$

Note that $h_n$ depends on $n$ as well as $C_1$, $C_2$, $\beta$ and $\gamma$. A similar definition applies in the setting of Definition 2.4.

REMARK 2.6. In the special case where $\gamma = 2$ and $\beta = 1$, the finite sample performance of this estimator has been investigated in Levine [25] together with the possible choice of bandwidth. A version of $K$-fold cross-validation has been recommended as the most suitable method. When utilized, it produces a variance estimator that in typical cases is not very sensitive to the choice of the mean function $g(x)$.

THEOREM 2.7. *Consider the nonparametric regression problem described by (1), with estimator as described in (8). Fix $C_1$, $C_2$, $\gamma > 0$ and $\beta > \gamma/(4\gamma + 2)$ to define functional classes $\mathcal{C}_\gamma$ and $\mathcal{C}_\beta$ according to the definition (2.2). Assume $p > \lfloor \gamma \rfloor$. Then the optimal bandwidth is $h_n \asymp n^{-1/(2\gamma+1)}$. Let $0 < \underline{a} \leq \overline{a} < \infty$. Then there are constants $\underline{B}$ and $\overline{B}$ such that*

(13) $$\underline{B} n^{-2\gamma/(2\gamma+1)} + o(n^{-2\gamma/(2\gamma+1)})$$
$$\leq R(V, \hat{V}) \leq \overline{B} n^{-2\gamma/(2\gamma+1)} + o(n^{-2\gamma/(2\gamma+1)})$$

*for all $h$ satisfying $\underline{a} \leq n^{1/(2\gamma+1)} h \leq \overline{a}$, uniformly for $g \in \mathcal{C}_\beta$, $V \in \mathcal{C}_\gamma$.*

Theorem 2.7 refers to properties of the integrated mean square error. Related results also hold for minimax risk at a point. The main results are stated in the following theorem.

THEOREM 2.8. *Consider the setting of Theorem 2.7. Let $x_0 \in (0, 1)$. Assume $p > \lfloor \gamma \rfloor$. Then the optimal bandwidth is $h_n(x) \asymp n^{-1/(2\gamma+1)}$. Let $0 < \underline{a} \leq \overline{a} < \infty$. Then there are constants $\underline{B}$ and $\overline{B}$ such that*

(14) $$\underline{B} n^{-2\gamma/(2\gamma+1)} + o(n^{-2\gamma/(2\gamma+1)}) \leq R(V(x_0), \hat{V}_{h_n}(x_0))$$
$$\leq \overline{B} n^{-2\gamma/(2\gamma+1)} + o(n^{-2\gamma/(2\gamma+1)})$$

for all $h(x)$ satisfying $\underline{a} \leq n^{1/(2\gamma+1)} h \leq \overline{a}$, uniformly for $g \in \mathcal{C}_\beta$, $V \in \mathcal{C}_\gamma$.

The proof of these theorems can be found in the Appendix. The minimax rates obtained in (13) and (14) will be shown in Theorems 4.1 and 4.2 to be optimal in the setting of Theorem 2.7. At this point, the following remarks may be helpful.

REMARK 2.9. If one assumes that $\beta = \gamma/(4\gamma + 2)$ in the definition of the functional class $\mathcal{C}_\beta$, the conclusions of Theorems 2.7 and 2.8 remain valid, but the constants $\underline{B}$ and $\overline{B}$ appearing in them become dependent on $\beta$. If $\beta < \gamma/(4\gamma + 2)$, the conclusion (14) does not hold. For more details, see comments preceding Theorem 4.2 and the Appendix.

REMARK 2.10. Müller and Stadtmüller [28] considered the general quadratic form based estimator similar to our (8) and derived convergence rates for its mean squared error. They also were the first to point out an error in the paper by Hall and Carroll [20] (see Müller and Stadtmüller [28], pages 214 and 221). They use a slightly different (more restrictive) definition of the classes $\mathcal{C}_\gamma$ and $\mathcal{C}_\beta$ and only establish rates of convergence and error terms on those rates for fixed functions $V$ and $g$ within the classes $\mathcal{C}_\gamma$ and $\mathcal{C}_\beta$. Our results resemble these but we also establish the rates of convergence *uniformly* over the functional classes $\mathcal{C}_\beta$ and $\mathcal{C}_\gamma$ and therefore our bounds are of the minimax type.

REMARK 2.11. It is important to notice that the asymptotic mean squared risks in Theorems 2.7 and 2.8 can be further reduced by proper choice of the difference sequence $\{d_j\}$. The proof in the Appendix supplemented with material in Hall, Kay and Titterington [21] shows that the asymptotic variance of our estimators will be affected by the choice of the difference sequence, but the choice of this sequence does not affect the bias in asymptotic calculations. The effect on the asymptotic variance is to multiply it by a constant proportional to

$$(15) \qquad C = 2\left(1 + 2 \sum_{k=1}^{r} \left( \sum_{j=0}^{r-1-k} d_j d_{j+k} \right)^2 \right).$$

For any given value of $r$ there is a difference sequence that minimizes this constant. A computational algorithm for these sequences is given in Hall, Kay and Titterington [21]. The resulting minimal constant as a function of $r$ is $C_{\min} = (2r+1)/r$.

**3. Asymptotic normality.** As a next step, we establish that the estimator (8) is asymptotically normal. We recall that the local polynomial regression estimator $\hat{V}_h(x)$ can be represented as

$$\hat{V}_h(x) = \sum_{i=\lfloor r/2 \rfloor + 1}^{n+\lfloor r/2 \rfloor - r} K_{n;h,x}(x_i) \Delta_{r,i}^2, \tag{16}$$

where $K_{n;h,x}(x_i) = K_{n,x}(\frac{x-x_i}{h})$. Here $K_{n,x}(\frac{x-x_i}{h})$ can be thought of as a centered and rescaled nonnegative local kernel function whose shape depends on the location of design points $x_i$, the point of estimation $x$ and the number of observations $n$. We know that $K_{n,x}(\frac{x-x_i}{h})$ satisfies discrete moment conditions,

$$\sum_{i=\lfloor r/2 \rfloor + 1}^{n+\lfloor r/2 \rfloor - r} K_{n,x}\left(\frac{x-x_i}{h}\right) = 1, \tag{17}$$

$$\sum_{i=\lfloor r/2 \rfloor + 1}^{n+\lfloor r/2 \rfloor - r} (x-x_i)^q K_{n,x}\left(\frac{x-x_i}{h}\right) = 0 \tag{18}$$

for any $q = 1, \ldots, p$. We also need the fact that the support of $K^n(\cdot)$ is contained in that of $K(\cdot)$; in other words, $K^n(\cdot) = 0$ whenever $|x_i - x| > h$. For more details see, for example, Fan and Gijbels [17]. Now we can state the following result.

THEOREM 3.1. *Consider the nonparametric regression problem described by (1), with estimator as described in (8). We assume that the functions $g(x)$ and $V(x)$ are continuous for any $x \in [0,1]$ and $V$ is bounded away from zero. Assume $\mu_{4+\nu} = E(\varepsilon_i)^{4+\nu} < \infty$ for some $\nu > 0$. Then, as $h \to 0$, $n \to \infty$ and $nh \to \infty$, we find that*

$$\sqrt{nh}(\hat{V}_h(x) - V(x) - O(h^{2\gamma})) \tag{19}$$

*is asymptotically normal with mean zero and variance $\sigma^2$ where $0 < \sigma^2 < \infty$.*

PROOF. To prove this result, we rely on the CLT for partial sums of a generalized linear process

$$X_n = \sum_{i=1}^{n} a_{ni} \xi_i, \tag{20}$$

where $\xi_i$ is a mixing sequence. This and several similar results were established in Peligrad and Utev [32]. Thus, the estimator (8) can be easily represented in the form (20) with $K_{n;h,x}(x_i)$ as $a_{ni}$. What remains is to verify the conditions of Theorem 2.2(c) in Peligrad and Utev [32].

- The first condition is

$$\max_{1 \leq i \leq n} |a_{ni}| \to 0 \tag{21}$$

as $n \to \infty$ and it is immediately satisfied since

$$K_{n;h,x}(x_i) = O((nh)^{-1}) \tag{22}$$

uniformly for all $x \in [0,1]$.

- The second condition is

$$\sup_n \sum_{i=1}^n a_{ni}^2 < \infty. \tag{23}$$

It can be verified by using the Cauchy–Schwarz inequality and (22).

- To establish uniform integrability of $\xi_i^2 \equiv \Delta_{r,i}^4$, we use a simple criterion mentioned in Shiryaev [37] that requires existence of the nonnegative, monotonically increasing function $G(t)$, defined for $t \geq 0$, such that

$$\lim_{t \to \infty} \frac{G(t)}{t} = \infty$$

and

$$\sup_i E[G(\Delta_{r,i}^4)] < \infty.$$

It is enough to choose $G(t) = t^\nu$ for small $\nu >$ to have this condition satisfied. Finally, the remaining three conditions of Peligrad and Utev [32] are trivially satisfied. □

**4. Asymptotic minimaxity and related issues.** Lower bounds on the asymptotic minimax rate for estimating a nonparametric variance in formulations related to that in (1) have occasionally been studied in earlier literature. Two papers seem particularly relevant. Munk and Ruymgaart [30] study a different, but related problem. Their paper contains a lower bound on the asymptotic minimax risk for their setting. In particular, their setting involves a problem with random design, rather than the fixed design case in (1). Their proof uses the Van Trees inequality and relies heavily on the fact that their $(X_i, Y_i)$ pairs are independent and identically distributed. While it may well be possible to do so, it is not immediately evident how to modify their argument to apply to the setting (1).

Hall and Carroll [20] consider a setting similar to ours. Their equation (2.13) claims (in our notation) that there is a constant $K > 0$, possibly depending on $C_1$, $C_2$, $\beta$ such that for any estimator $\tilde{V}$

$$\sup\{R(V(x_0), \tilde{V}(x_0)) : V \in \mathcal{C}_\gamma, g \in \mathcal{C}_\beta\} \geq K \max\{n^{-2\gamma/(2\gamma+1)}, n^{-4\beta/(2\beta+1)}\}. \tag{24}$$

Note that $n^{-2\gamma/(2\gamma+1)} = o(n^{-4\beta/(2\beta+1)})$ for $\beta < \gamma/(2\gamma + 2)$. It thus follows from (14) in our Theorem 2.8 that for any $\gamma/(4\gamma + 2) < \beta < \gamma/(2\gamma + 2)$ and $n$ sufficiently large

$$\begin{aligned}(25)\quad &\sup\{R(V(x_0), \hat{V}_{h_n}(x_0)) : V \in \mathcal{C}_\gamma, g \in \mathcal{C}_\beta\} \\ &\ll K \max\{n^{-2\gamma/(2\gamma+1)}, n^{-4\beta/(2\beta+1)}\},\end{aligned}$$

where $h_n$ is yet again the optimal bandwidth. This contradicts the assertion in Hall and Carroll [20], and shows that their assertion (2.13) is in error— as is the argument supporting it that follows (C.3) of their article. For a similar commentary see also Müller and Stadtmüller [28]. Because of this contradiction it is necessary to give an independent statement and proof of a lower bound for the minimax risk. That is the goal of this section, where we treat the case in which $\beta \geq \gamma/(4\gamma + 2)$. The minimax lower bound for the case in which $\beta < \gamma/(4\gamma + 2)$ requires different methods which are more sophisticated. That case, as well as some further generalizations, have been treated in Wang, Brown, Cai and Levine [42] as a sequel to the present paper. That paper proves ratewise sharp lower and upper bounds for the case where $\beta < \gamma/(4\gamma + 2)$.

We have treated both mean squared error at a point (in Theorem 2.8) and integrated mean squared error (in Theorem 2.7). Correspondingly, we provide statements of lower bounds on the minimax rate for each of these cases. The local version of the lower bound result for the minimax risk is obtained under the assumption of normality of errors $\varepsilon_i$. See Section 2 for the definition of $R$ and other quantities that appear in the following statements.

THEOREM 4.1. *Consider the nonparametric regression problem described by (1). Fix $C_1$, $C_2$, $\beta$ and $\gamma$ to define functional classes $\mathcal{C}_\gamma$, $\mathcal{C}_\beta$ according to (2.2). Also assume that $\varepsilon_i \sim N(0,1)$ and independent. Then there is a constant $K > 0$ such that*

$$(26) \quad \inf\{\sup\{R(V, \tilde{V}) : V \in \mathcal{C}_\gamma^+, g \in \mathcal{C}_\beta\} : \tilde{V}\} \geq K n^{-2\gamma/(2\gamma+1)}$$

*where the* inf *is taken over all possible estimators of the variance function $V$.*

Our argument relies on the so-called "two-point" argument, introduced and extensively analyzed in Donoho and Liu [11, 12].

THEOREM 4.2. *Consider the nonparametric regression problem described by (1). Fix $C_1$, $C_2$, $\beta$ and $\gamma$ to define functional classes $\mathcal{C}_\gamma$, $\mathcal{C}_\beta$ according to (2.2). Also assume that $\varepsilon_i \sim N(0,1)$ and independent. Then there is a constant $K > 0$ such that*

$$(27) \quad \inf\{\sup\{R(V(x_0)), \tilde{V}(x_0)) : V \in \mathcal{C}_\gamma, g \in \mathcal{C}_\beta\} : \tilde{V}\} \geq K n^{-2\gamma/(2\gamma+1)}$$

*where the* inf *is taken over all possible estimators of the variance function $V$.*

PROOF. It is easier to begin with the proof of Theorem 4.2 and then proceed to the proof of Theorem 4.1. We will use a two-point modulus-of-continuity argument to establish the lower bound. Such an argument was pioneered by Donoho and Liu [11, 12] for a different though related problem. See also Hall and Carroll [20] and Fan [16].

We assume without loss of generality that $g \equiv 0$. Define the function

$$h(t) = \begin{cases} 2 - |t|^\gamma, & \text{if } 0 \le |t| \le 1, \\ (2 - |t|)^\gamma, & \text{if } 1 < |t| \le 2, \\ 0, & \text{if } |t| > 2. \end{cases} \tag{28}$$

Assume (for convenience only) that $C_1 > 2$. Let $d$ be a constant satisfying $0 < d < C_2$ and let

$$f_{\delta,l}(x) = d + l\delta h\left(\frac{x - x_0}{\delta^{1/\gamma}}\right). \tag{29}$$

Then $f_{\varepsilon,\pm 1} \in \mathcal{C}_\gamma$ for $\delta > 0$ sufficiently small. Let $H$ denote the Hellinger distance between densities, that is, for any two probability densities $m_1$, $m_2$ dominated by a measure $\mu(dz)$,

$$H^2(m_1, m_2) = \int (\sqrt{m_1(z)} - \sqrt{m_2(z)})^2 \mu(dz). \tag{30}$$

Here are two basic facts about this metric that will be used below. If $Z = \{Z_j : j - 1, \ldots, n\}$ where the $Z_j$ are independent with densities $\{m_{kj} : j = 1, \ldots, n\}$, $k = 1, 2$ and $m_k = \Pi_j m_{kj}$ denotes the product density, then

$$H^2(m_1, m_2) \le \sum_j H^2(m_{1j}, m_{2j}); \tag{31}$$

and if $m_i$ are univariate normal densities with mean 0 and variance $\sigma_i^2$, $i = 1, 2$, then

$$H^2(m_1, m_2) \le 2\left(\frac{\sigma_1^2}{\sigma_2^2} - 1\right)^2. \tag{32}$$

For more details see Brown and Low [3] and Brown et al. [1].

It follows that if $m_k$, $k = 1, 2$, are the joint densities of the observations $\{x_i, Y_i, i = 1, \ldots, n\}$ of (1) with $g \equiv 0$ and $f_k = f_{\delta,(-1)^k}$ then

$$\begin{aligned} H^2(m_1, m_2) &\le \sum_i 2\left(\frac{f_{\delta,-1}(x_i)}{f_{\delta,1}(x_i)} - 1\right)^2 \\ &\le 8 \sum_i \delta^2 h^2\left(\frac{x_i - x_0}{\delta^{1/\gamma}}\right) = O(n\gamma^{(2\gamma+1)/\gamma}). \end{aligned} \tag{33}$$

For this setting the Hellinger modulus-of-continuity, $\omega(\cdot)$ (Donoho and Liu [12], equation (1.1)), is defined as the inverse function corresponding to the value $H(m_1, m_2)$. Hence it satisfies

$$\omega^{-1}(\gamma) = O(n^{1/2} \gamma^{(2\gamma+1)/2\gamma}). \tag{34}$$

Equation (27) then follows, as established in Donoho and Liu [12]. Although this completes the proof of Theorem 4.2, we also provide a sketch of the argument based on (34). See Donoho and Liu [12] and references cited therein for more details. □

PROOF OF THEOREM 4.1. We omit this proof for the sake of brevity. It begins from the result in Theorem 4.2 and then follows along the lines first described in detail in Donoho, Liu and MacGibbon [13]. This theorem can be also viewed as a consequence of the general results on the global convergence of nonparametric estimators by Stone [39] and Efromovich [14] that do not require normality of errors $\varepsilon_i$. □

## APPENDIX

PROOFS OF THEOREMS 2.7 AND 2.8. Fix $r$ and functional classes $\mathcal{C}_\gamma$ and $\mathcal{C}_\beta$. For the sake of brevity, we write $\Delta_i \equiv \Delta_{r,i}$. Our main tools in this proof are the representation (16) of the variance estimator $\hat{V}_h(x)$ and the properties (17)–(18). We also use the property

$$\sum_{i=\lfloor r/2 \rfloor - 1}^{n+\lfloor r/2 \rfloor - r} (K_{n;h,x}(x_i))^2 = O\left(\frac{1}{nh}\right). \tag{35}$$

(35) follows from (22) and the Cauchy–Schwarz inequality. Here and later, $O$ is uniform for all $V \in \mathcal{C}_\gamma$, $g \in \mathcal{C}_\beta$ and $\{h\} = \{h_n\}$. Now,

$$E(\Delta_i^2) = \text{Var}(\Delta_i) + (E(\Delta_i))^2, \tag{36}$$

where

$$\text{Var}(\Delta_i) = \sum d_j^2 \text{Var}(y_{j+i-\lfloor r/2 \rfloor}) = V(x_i) + O\left(\left(\frac{1}{n}\right)^\gamma\right) \tag{37}$$

and

$$E(\Delta_i) = O\left(\left(\frac{1}{n}\right)^\beta\right) \tag{38}$$

since $\sum d_j = 0$, $\sum d_j^2 = 1$ and $x_{i+r-\lfloor r/2 \rfloor} - x_{i-\lfloor r/2 \rfloor} = O(\frac{1}{n})$. This provides an asymptotic bound on the bias as

$$\begin{aligned}
\text{Bias } \hat{V}_h(x) &= \sum_{i=\lfloor r/2 \rfloor + 1}^{n+\lfloor r/2 \rfloor - r} (V(x_i) - V(x)) K_{n;h,x}(x_i) + O(n^{-\gamma}) + O(n^{-\beta}) \\
&= O(h^\gamma) + O(n^{-\gamma}) + O(n^{-\beta}).
\end{aligned} \tag{39}$$

The last step in (39) is a very minor variation of the technique employed in Wang, Brown, Cai and Levine [42] (see pages 10–11).

Next, we need to use the fact that $\Delta_i$ and $\Delta_j$ are independent if $|i - j| \geq r + 1$. Hence,

$$\operatorname{Var} \hat{V}_h(x) = \operatorname{Var}\left(\sum_{i=\lfloor r/2 \rfloor+1}^{n+\lfloor r/2 \rfloor-r} K_{n;h,x}(x_i)\Delta_i^2\right)$$

$$= \sum_{i=\lfloor r/2 \rfloor+1}^{n+\lfloor r/2 \rfloor-r} \sum_{j=i-r}^{i+r} K_{n;h,x}(x_i)K_{n;h,x}(x_j) \operatorname{Cov}(\Delta_i^2, \Delta_j^2)$$

$$\leq \sum_{i=\lfloor r/2 \rfloor+1}^{n+\lfloor r/2 \rfloor-r} \sum_{j=i-r}^{i+r} 4^{-1}((K_{n;h,x}(x_i))^2 + (K_{n;h,x}(x_j))^2)$$

$$\times (\operatorname{Var}\Delta_i^2 + \operatorname{Var}\Delta_j^2)$$

It is easy to see that

$$\Delta_i^2 = \left(\sum_{j=0}^{r} d_j y_{j+i-\lfloor r/2 \rfloor}\right)^2$$

$$= \left(\sum_{j=0}^{r} d_j \sqrt{V(x_{j+i-\lfloor r/2 \rfloor})}\varepsilon_{i+j-\lfloor r/2 \rfloor} + O(n^{-\beta})\right)^2,$$

and this means, in turn, that

$$\operatorname{Var}\Delta_i^2 \leq C_2^2 \operatorname{Var}\left(\sum_{j=0}^{r} d_j \varepsilon_{i+j-\lfloor r/2 \rfloor} + O(n^{-\beta})\right)^2$$

$$\leq C_2^2(r+1)\mu_4 + O(n^{-2\beta}) + O(n^{-4\beta}) = O(1).$$

Hence,

$$\operatorname{Var} \hat{V}_h(x) \leq O(1) \sum_{i=\lfloor r/2 \rfloor+1}^{n+\lfloor r/2 \rfloor-r} \sum_{j=i-r}^{i+r} ((K_{n;h,x}(x_i))^2 + (K_{n;h,x}(x_j))^2)$$

(40)
$$= O\left(\frac{1}{nh}\right).$$

Combining the bounds in (39) and (40) yields the assertion of the theorem since $2\beta > \gamma/(2\gamma + 1)$. □

**Acknowledgments.** We wish to thank T. Cai and L. Wang for pointing out the significance of the article of Hall and Carroll [20] and its relation to our (14).

Department of Statistics  
The Wharton School  
University of Pennsylvania  
Philadelphia, Pennsylvania 19104-6340  
USA  
E-mail: lbrown@wharton.upenn.edu

Department of Statistics  
Purdue University  
150 N. University Street  
West Lafayette, Indiana 47907  
USA  
E-mail: mlevins@stat.purdue.edu